\input amstex
\documentstyle{amsppt}
%----------------------------------------------------------------
% Title:     A rough classification of potentially invertible cubic
%            transformations of tye real plane 
%            of the real plane.
% Author:    Ruslan Sharipov
% Comments:  AmSTeX, 8 pages, amsppt style
% MSC-class: 14E05, 26B10, 57S25
%----------------------------------------------------------------
%           Replacement for output macro definition
%
\catcode`@=11
\redefine\output@{%
  \def\break{\penalty-\@M}\let\par\endgraf
  \ifodd\pageno\global\hoffset=105pt\else\global\hoffset=8pt\fi  
  \shipout\vbox{%
    \ifplain@
      \let\makeheadline\relax \let\makefootline\relax
    \else
      \iffirstpage@ \global\firstpage@false
        \let\rightheadline\frheadline
        \let\leftheadline\flheadline
      \else
        \ifrunheads@ %\let\makefootline\relax
        \else \let\makeheadline\relax
        \fi
      \fi
    \fi
    \makeheadline \pagebody \makefootline}%
  \advancepageno \ifnum\outputpenalty>-\@MM\else\dosupereject\fi
}
\def\Beta{\mathchar"0\hexnumber@\rmfam 42}
\catcode`\@=\active
%----------------------------------------------------------------
\nopagenumbers
\chardef\textvolna='176

\chardef\bigalpha='013
\def\negskp{\hskip -2pt}

\chardef\degree="5E
\def\compos{\,\raise 1pt\hbox{$\sssize\circ$} \,}

\def\RP{\operatorname{\Bbb R\text{\rm P}}}
%\font\eightrm=cmr8
%\def\LT{\operatorname{\text{\eightrm LT}}}
%\def\LM{\operatorname{\text{\eightrm LM}}}
%\def\LC{\operatorname{\text{\eightrm LC}}}
%\accentedsymbol\hatgamma{\kern 2pt\hat{\kern -2pt\gamma}}
%\accentedsymbol\checkgamma{\kern 2.5pt\check{\kern -2.5pt\gamma}}
\def\blue#1{#1}

\catcode`#=11\def\diez{#}\catcode`#=6
\catcode`&=11\catcode`&=4
\catcode`_=11\def\podcherkivanie{_}\catcode`_=8
\catcode`\^=11\catcode`\^=7
\catcode`~=11\catcode`~=\active
\def\mycite#1{\cite{\blue{#1}}\immediate\special{ps:
     ShrHPSdict begin /ShrBORDERthickness 0 def}}
\def\myciterange#1#2#3#4{\cite{\blue{#2#3#4}}\immediate\special{ps:
     ShrHPSdict begin /ShrBORDERthickness 0 def}}
\def\mytag#1{%
    \tag#1}
\def\mythetag#1{\thetag{\blue{#1}}\immediate\special{ps:
     ShrHPSdict begin /ShrBORDERthickness 0 def}}
\def\myrefno#1{\no#1}
\def\myhref#1#2{\blue{#2}\immediate\special{ps:
     ShrHPSdict begin /ShrBORDERthickness 0 def}}
\def\myEarXivlink{\myhref{http://arXiv.org}{http:/\negskp/arXiv.org}}

\def\mytheorem#1{\csname proclaim\endcsname{Theorem #1}}
\def\mytheoremwithtitle#1#2{\csname proclaim\endcsname{Theorem #1#2}}
\def\mythetheorem#1{\blue{#1}\immediate\special{ps:
     ShrHPSdict begin /ShrBORDERthickness 0 def}}
\def\mylemma#1{\csname proclaim\endcsname{Lemma #1}}
\def\mylemmawithtitle#1#2{\csname proclaim\endcsname{Lemma #1#2}}

\def\mycorollary#1{\csname proclaim\endcsname{Corollary #1}}

\def\mydefinition#1{\definition{Definition #1}}
\def\mythedefinition#1{\blue{#1}\immediate\special{ps:
     ShrHPSdict begin /ShrBORDERthickness 0 def}}
\def\myconjecture#1{\csname proclaim\endcsname{Conjecture #1}}
\def\myconjecturewithtitle#1#2{\csname proclaim\endcsname{Conjecture #1#2}}

\def\myproblem#1{\csname proclaim\endcsname{Problem #1}}
\def\myproblemwithtitle#1#2{\csname proclaim\endcsname{Problem #1#2}}

%----------------------------------------------------------------
% Cyrillic fonts definition

%\font\tencyr=wncyr10
%----------------------------------------------------------------
\pagewidth{360pt}
\pageheight{606pt}
\topmatter
\title
A rough classification of potentially invertible cubic transformations
of the real plane.
\endtitle
\rightheadtext{A rough classification \dots}
\author
Ruslan Sharipov
\endauthor
\address Bashkir State University, 32 Zaki Validi street, 450074 Ufa, Russia
\endaddress
\email
\myhref{mailto:r-sharipov\@mail.ru}{r-sharipov\@mail.ru}
\endemail
\abstract
     A polynomial transformation of the real plane $\Bbb R^2$ is a mapping
$\Bbb R^2\to\Bbb R^2$ given by two polynomials of two variables. Such a
transformation is called cubic if the degrees of its polynomials are not
greater than three. In the present paper a rough classification scheme for
cubic transformations of $\Bbb R^2$ is suggested. It is based on quartic 
forms associated with these transformations. 
\endabstract
\endtopmatter
%\loadbold
%\loadeufb
\TagsOnRight
\document
%\input countstyle

%\special{header=resource.eps}
\head
1. Introduction.
\endhead
     Let $f\!:\,\Bbb R^2\to\Bbb R^2$ be a cubic transformation of $\Bbb R^2$.
In the coordinate presentation it is given by the following formula:
$$
\hskip -2em
y^{\kern 0.5pt i}=\sum^2_{m=1}\sum^2_{n=1}\sum^2_{p=1}
F^{\kern 0.5pt i}_{m\kern 0.1pt np}\,x^m\,x^n\,x^{\kern 0.1pt p}
+\ldots.
\mytag{1.1}
$$
We use upper and lower indices in \mythetag{1.1} according to Einstein's 
tensorial notation (see \mycite{1}). Polynomial mappings $f\!:\,\Bbb R^2
\to\Bbb R^2$ were considered in \mycite{2}, \mycite{3} in connection with
the real Jacobian conjecture. The rational real Jacobian conjecture is its 
generalization (see \mycite{4}). Polynomial mappings $f\!:\,\Bbb C^2\to\Bbb 
C^2$ of the complex plane $\Bbb C^2$ were considered in 
\myciterange{1}{5}{--}{9}.\par
     My interest to cubic polynomial transformations of $\Bbb R^2$ is due to 
their application to the perfect cuboid problem (see \mycite{10} and 
\mycite{11}). However, in the present paper we study these transformations 
by themselves. Cubic polynomial transformations of the form \mythetag{1.1}
are associated with some definite quartic forms (see \mycite{12}). They are 
produced using components of the tensor $F$ in \mythetag{1.1}. In the present
paper we use them as a basis for building a rough classification of
potentially invertible cubic transformations of $\Bbb R^2$. This classification
can further be refined using other invariants and/or methods. The classification
of potentially invertible quadratic transformations of $\Bbb R^2$ can be found
in \mycite{13}.\par 
\head
2. Determinants and quartic forms. 
\endhead
     The tensor $F$ determines the leading terms of the cubic transformation 
\mythetag{1.1}. They are of the most attention below. Terms of lower degrees 
are denoted by dots. They are not used in our rough classification scheme.\par 
     Along with \mythetag{1.1}, we consider linear transformations 
$\varphi\!:\,\Bbb R^2\to\Bbb R^2$ of the form
$$
\hskip -2em
y^{\kern 0.5pt i}=\sum^2_{m=1}T^{\kern 0.5pt i}_m\,x^m+a^i.
\mytag{2.1}
$$
\mydefinition{1.1} Two cubic mappings $f\!:\,\Bbb R^2\to \Bbb R^2$ and
$\tilde f\!:\,\Bbb R^2\to \Bbb R^2$ are called equivalent if there are two 
invertible linear transformations $\varphi_1\!:\,\Bbb R^2\to \Bbb R^2$ and
$\varphi_2\!:\,\Bbb R^2\to \Bbb R^2$ such that $\varphi_1\compos f=
\tilde f\compos\varphi_2$.
\enddefinition
     It is natural to classify cubic transformations up to the
equivalence introduced in Definition~\mythedefinition{1.1}. For this
purpose we need some parameters which are relatively stable when 
passing from a given cubic transformations to an equivalent one. These
parameters are constructed through the tensor $F$ in \mythetag{1.1}.
\par
     The tensor $F$ in \mythetag{1.1} is symmetric with respect to its lower
indices. Taking into account this symmetry we can write \mythetag{1.1} as 
$$
\aligned
&y^1=F^{\kern 1pt 1}_{111}\,(x^1)^3+3\,F^{\kern 1pt 1}_{112}\,(x^1)^2\,x^2
+3\,F^{\kern 1pt 1}_{122}\,x^1\,(x^2)^2+F^{\kern 1pt 1}_{222}\,(x^2)^3+\dots,\\
&y^2=F^{\kern 1pt 2}_{111}\,(x^1)^3+3\,F^{\kern 1pt 2}_{112}\,(x^1)^2\,x^2
+3\,F^{\kern 1pt 2}_{122}\,x^1\,(x^2)^2+F^{\kern 1pt 2}_{222}\,(x^2)^3+\dots.
\endaligned
\mytag{2.2}
$$
According to \mycite{12} and \mycite{13}, we consider the following 
determinants:
$$
\xalignat 2
&\hskip -2em
G_{1111}=\vmatrix F^{\kern 1pt 1}_{111} & F^{\kern 1pt 1}_{112}\\
\vspace{2ex}
F^{\kern 1pt 2}_{111} & F^{\kern 1pt 2}_{112}\endvmatrix,
&&G_{1112}=\vmatrix F^{\kern 1pt 1}_{111} & F^{\kern 1pt 1}_{122}\\
\vspace{2ex}
F^{\kern 1pt 2}_{111} & F^{\kern 1pt 2}_{122}\endvmatrix,\\
\vspace{2ex}
&\hskip -2em
G_{1122}=\vmatrix F^{\kern 1pt 1}_{111} & F^{\kern 1pt 1}_{222}\\
\vspace{2ex}
F^{\kern 1pt 2}_{111} & F^{\kern 1pt 2}_{222}\endvmatrix,
&&G_{1212}=\vmatrix F^{\kern 1pt 1}_{112} & F^{\kern 1pt 1}_{122}\\
\vspace{2ex}
F^{\kern 1pt 2}_{112} & F^{\kern 1pt 2}_{122}\endvmatrix,
\mytag{2.3}\\
\vspace{2ex}
&\hskip -2em
G_{1222}=\vmatrix F^{\kern 1pt 1}_{112} & F^{\kern 1pt 1}_{222}\\
\vspace{2ex}
F^{\kern 1pt 2}_{112} & F^{\kern 1pt 2}_{222}\endvmatrix,
&&G_{2222}=\vmatrix F^{\kern 1pt 1}_{122} & F^{\kern 1pt 1}_{222}\\
\vspace{2ex}
F^{\kern 1pt 2}_{122} & F^{\kern 1pt 2}_{222}\endvmatrix.
\endxalignat
$$
The determinants \mythetag{2.3} are related to six quartic forms 
$\omega[1]$, $\omega[2]$, $\omega[3]$, $\omega[4]$, $\omega[5]$, 
$\omega[6]$ introduced in \mycite{12}. They are given by the following 
formulas:
$$
\allowdisplaybreaks
\gather
\hskip -2em
\gathered
\omega[1]=G_{1111}\,(z^1)^4+2\,G_{1112}\,(z^1)^3\,z^2\,+\\
+\,(3\,G_{1212}+G_{1122})\,(z^1)^2\,(z^2)^2+2\,G_{1222}\,z^1\,(z^2)^3
+G_{2222}\,(z^2)^4,
\endgathered
\mytag{2.4}\\
\vspace{2ex}
\hskip -2em
\gathered
\omega[2]=2\,G_{1111}\,(z^1)^3\,z^3+G_{1112}\,(z^1)^3\,z^4
+3\,G_{1112}\,(z^1)^2\,z^2\,z^3\,+\\
+\,(3\,G_{1212}+G_{1122})\,(z^1)^2\,z^2\,z^4
+(3\,G_{1212}+G_{1122})\,z^1\,(z^2)^2\,z^3\,+\,\\
+\,3\,G_{1222}\,z^1\,(z^2)^2\,z^4
+G_{2222}\,\,(z^2)^3\,z^3+2\,G_{2222}\,\,(z^2)^3\,z^4,
\endgathered
\mytag{2.5}\\
\vspace{2ex}
\hskip -2em
\gathered
\omega[3]=3\,G_{1111}\,(z^1)^2\,(z^3)^2+3\,G_{1112}\,(z^1)^2\,z^3\,z^4
+G_{1122}\,(z^1)^2\,(z^4)^2\,+\\
+\,3\,G_{1112}\,z^1\,z^2\,(z^3)^2
+(9\,G_{1212}+G_{1122})\,z^1\,z^2\,z^3\,z^4+3\,G_{1222}\,z^1\,z^2\,(z^4)^2\,+\\
+\,G_{1122}\,(z^2)^2\,(z^3)^2
+3\,G_{1222}\,\,(z^2)^2\,z^3\,z^4+3\,G_{2222}\,\,(z^2)^2\,(z^4)^2,
\endgathered
\mytag{2.6}\\
\hskip -2em
\gathered
\omega[4]=G_{1111}\,(z^1)^2\,(z^3)^2+G_{1112}\,(z^1)^2\,z^3\,z^4
+G_{1212}\,(z^1)^2\,(z^4)^2\,+\\
+\,G_{1112}\,z^1\,z^2\,(z^3)^2
+\,(G_{1212}+G_{1122})\,z^1\,z^2\,z^3\,z^4
+G_{1222}\,z^1\,z^2\,(z^4)^2\,+\\
+\,G_{1212}\,(z^2)^2\,(z^3)^2+G_{1222}\,(z^2)^2\,z^3\,z^4 
+G_{1222}\,(z^2)^2\,(z^4)^2,
\endgathered
\mytag{2.7}\\
\vspace{2ex}
\hskip -2em
\gathered
\omega[5]=2\,G_{1111}\,z^1\,(z^3)^3+G_{1112}\,z^2\,(z^3)^3
+3\,G_{1112}\,z^1\,(z^3)^2\,z^4\,+\\
+\,(3\,G_{1212}+G_{1122})\,z^2\,(z^3)^2\,z^4
+(3\,G_{1212}+G_{1122})\,z^1\,z^3\,(z^4)^2\,+\\
+\,3\,G_{1222}\,z^2\,z^3\,(z^4)^2
+G_{1222}\,z^1\,(z^4)^3
+2\,G_{2222}\,z^2\,(z^4)^3,
\endgathered
\mytag{2.8}\\
\vspace{2ex}
\hskip -2em
\gathered
\omega[6]=G_{1111}\,(z^3)^4+2\,G_{1112}\,(z^3)^3\,z^4\,+\\
+\,(3\,G_{1212}+G_{1122})\,(z^3)^2\,(z^4)^2
+2\,G_{1222}\,z^3\,(z^4)^3+G_{2222}\,(z^4)^4.
\endgathered
\mytag{2.9}
\endgather
$$
The determinants \mythetag{2.3} arise not only as coefficients in
the forms \mythetag{2.4}, \mythetag{2.5}, \mythetag{2.6}, \mythetag{2.7}, 
\mythetag{2.8}, \mythetag{2.9}. They generate these forms 
according to the following theorem.
\mytheorem{2.1} If a cubic transformation $f$ is produced as the
right composition $f=\tilde f\compos\varphi$ of another cubic 
transformation $\tilde f$ with a liner transformation $\varphi$ of 
the form \mythetag{2.1}, then the determinants \mythetag{2.3} associated 
with $f$ are expressed through the values of the quartic forms 
$\tilde\omega[1]$, $\tilde\omega[2]$, $\tilde\omega[3]$, $\tilde\omega[4]$, 
$\tilde\omega[5]$, $\tilde\omega[6]$ associated with the second cubic
transformation $\tilde f$ according to the formulas 
$$
\aligned
&G_{1111}(f)=\det T\cdot\tilde\omega[1](z^1,z^2),\\
&G_{1112}(f)=\det T\cdot\tilde\omega[2](z^1,z^2,z^3,z^4),\\
&G_{1122}(f)=\det T\cdot\tilde\omega[3](z^1,z^2,z^3,z^4),\\
&G_{1212}(f)=\det T\cdot\tilde\omega[4](z^1,z^2,z^3,z^4),\\
&G_{1222}(f)=\det T\cdot\tilde\omega[5](z^1,z^2,z^3,z^4),\\
&G_{2222}(f)=\det T\cdot\tilde\omega[6](z^3,z^4),
\endaligned
$$
where $z^1$, $z^2$, $z^3$, $z^4$ are given by the components of the matrix 
$T$ in \mythetag{2.1}:
$$
\xalignat 2
&\hskip -2em
z^1=T^1_1,
&&z^3=T^1_2,\\
\vspace{-1.5ex}
\mytag{2.10}\\
\vspace{-1.5ex}
&\hskip -2em
z^2=T^2_1,
&&z^4=T^2_2.\\
\endxalignat
$$
\endproclaim 
     The transformation of the determinants \mythetag{2.3} under left 
compositions is more simple. It is described by another theorem. 
\mytheorem{2.2} If a cubic transformation $f$ is produced as the
left composition $f=\varphi^{-1}\compos\tilde f$ of another cubic 
transformation $\tilde f$ with the inverse of a liner transformation 
$\varphi$ in \mythetag{2.1}, then the associated determinants 
of the transformations $f$ and $\tilde f$ are related to each other
according to the formulas 
$$
\xalignat 2
&\hskip -2em
G_{1111}=\det S\cdot\tilde G_{1111}, 
&&G_{1112}=\det S\cdot\tilde G_{1112},\\
&\hskip -2em
G_{1122}=\det S\cdot\tilde G_{1122}, 
&&G_{1212}=\det S\cdot\tilde G_{1212},
\mytag{2.11}\\
&\hskip -2em
G_{1222}=\det S\cdot\tilde G_{1222}, 
&&G_{2222}=\det S\cdot\tilde G_{2222},
\endxalignat
$$
where $S=T^{-1}$ is the inverse of the matrix $T$ \pagebreak whose components 
are used in \mythetag{2.1}.
\endproclaim 
     Theorems~\mythetheorem{2.1} and \mythetheorem{2.2} are proved by means
of direct calculations. Theorem~\mythetheorem{2.1} was formulated in \mycite{12}. 
The formulas \mythetag{2.11} were also written in \mycite{12}, though they
were not formulated as a theorem.\par
\head
3. The case where $\omega[1]$ is zero.
\endhead
     If the quartic form $\omega[1]$ in \mythetag{2.4} is identically zero, this 
means that the determinants $G_{1111}$, $G_{1112}$, $G_{1222}$, $G_{2222}$
in \mythetag{2.3} are zero:
$$
\xalignat 2
&\hskip -2em
G_{1111}=0, 
&&G_{1112}=0,\\
\vspace{-1.5ex}
\mytag{3.1}\\
\vspace{-1.5ex}
&\hskip -2em
G_{1222}=0, 
&&G_{2222}=0.
\endxalignat
$$
Apart from \mythetag{3.1}, this means that the following equality is fulfilled:
$$
\hskip -2em
3\,G_{1212}+G_{1122}=0. 
\mytag{3.2}
$$
The matrices producing the determinants $G_{1111}$, $G_{1112}$,
$G_{1122}$, $G_{1212}$, $G_{1222}$, $G_{2222}$ in \mythetag{2.3} 
are composed with the use of the following four vector-columns:
$$
\xalignat 4
&\hskip -2em
\Vmatrix F^{\kern 1pt 1}_{111}\\ 
\vspace{2ex}F^{\kern 1pt 2}_{111}\endVmatrix,
&&\Vmatrix F^{\kern 1pt 1}_{112}\\ 
\vspace{2ex} F^{\kern 1pt 2}_{112}\endVmatrix,
&&\Vmatrix F^{\kern 1pt 1}_{122}\\ 
\vspace{2ex} F^{\kern 1pt 2}_{122}\endVmatrix,
&&\Vmatrix F^{\kern 1pt 1}_{222}\\ 
\vspace{2ex} F^{\kern 1pt 2}_{222}\endVmatrix.
\quad
\mytag{3.3}
\endxalignat
$$
At least one of the vector-columns \mythetag{3.3} is nonzero (since
$f$ is cubic).\par 
    Assume that the first column in \mythetag{3.3} is nonzero, i\.\,e\.
assume that 
$$
\hskip -2em
\Vmatrix F^{\kern 1pt 1}_{111}\\ 
\vspace{2ex}F^{\kern 1pt 2}_{111}\endVmatrix\neq 0.
\mytag{3.4}
$$
Then due to $G_{1111}=0$ and $G_{1112}=0$ in \mythetag{3.1} the
second and the third columns in \mythetag{3.3} are expressed through 
the first one as follows:
$$
\xalignat 2
&\hskip -2em
\Vmatrix F^{\kern 1pt 1}_{112}\\ 
\vspace{2ex} F^{\kern 1pt 2}_{112}\endVmatrix
=\alpha\Vmatrix F^{\kern 1pt 1}_{111}\\ 
\vspace{2ex}F^{\kern 1pt 2}_{111}\endVmatrix,
&&\Vmatrix F^{\kern 1pt 1}_{122}\\ 
\vspace{2ex} F^{\kern 1pt 2}_{122}\endVmatrix
=\beta\Vmatrix F^{\kern 1pt 1}_{111}\\ 
\vspace{2ex}F^{\kern 1pt 2}_{111}\endVmatrix.
\mytag{3.5}
\endxalignat
$$
Substituting \mythetag{3.5} into the determinant $G_{1212}$ in 
\mythetag{2.3} we derive $G_{1212}=0$. Applying $G_{1212}=0$ to
\mythetag{3.2}, we derive the equality 
$$
\hskip -2em
G_{1122}=0.
\mytag{3.6}
$$
The equality \mythetag{3.6} combined with the inequality \mythetag{3.4} 
means that the fourth column in \mythetag{3.3} is expressed through the 
first one as follows:
$$
\hskip -2em
\Vmatrix F^{\kern 1pt 1}_{222}\\ 
\vspace{2ex} F^{\kern 1pt 2}_{222}\endVmatrix
=\gamma\Vmatrix F^{\kern 1pt 1}_{111}\\ 
\vspace{2ex}F^{\kern 1pt 2}_{111}\endVmatrix. 
\mytag{3.7}
$$\par
     The equalities \mythetag{3.5} and \mythetag{3.7} mean that the cubic 
terms in \mythetag{2.2} are proportional to each other. Therefore, applying
the left composition with some properly chosen linear transformation of the 
form \mythetag{2.1}, we can pass to an equivalent cubic transformation 
such that $F^{\kern 1pt 2}_{111}=F^{\kern 1pt 2}_{112}=F^{\kern 1pt 2}_{122}
=F^{\kern 1pt 2}_{222}=0$.\par
      If the inequality \mythetag{3.4} is not fulfilled, then at least one
of the last three columns in \mythetag{3.3} is nonzero. Assume that the
second column in \mythetag{3.3} is nonzero: 
$$
\xalignat 2
&\hskip -2em
\hskip -2em
\Vmatrix F^{\kern 1pt 1}_{111}\\ 
\vspace{2ex}F^{\kern 1pt 2}_{111}\endVmatrix=0,
&&\Vmatrix F^{\kern 1pt 1}_{112}\\ 
\vspace{2ex} F^{\kern 1pt 2}_{112}\endVmatrix\neq 0.
\mytag{3.8}
\endxalignat
$$
From \mythetag{3.8} we derive $G_{1122}=0$. Combining $G_{1122}=0$ with 
\mythetag{3.2}, we derive $G_{1212}=0$. The equality $G_{1212}=0$ combined
with $G_{1222}=0$ in \mythetag{3.1} and with \mythetag{3.8} means that the 
last two columns \mythetag{3.3} are expressed through the second one: 
$$
\xalignat 2
&\hskip -2em
\Vmatrix F^{\kern 1pt 1}_{122}\\ 
\vspace{2ex} F^{\kern 1pt 2}_{122}\endVmatrix
=\alpha\Vmatrix F^{\kern 1pt 1}_{112}\\ 
\vspace{2ex} F^{\kern 1pt 2}_{112}\endVmatrix,
&&\Vmatrix F^{\kern 1pt 1}_{222}\\ 
\vspace{2ex} F^{\kern 1pt 2}_{222}\endVmatrix
=\beta\Vmatrix F^{\kern 1pt 1}_{112}\\ 
\vspace{2ex} F^{\kern 1pt 2}_{112}\endVmatrix.
\mytag{3.9}
\endxalignat
$$
The relationships \mythetag{3.8} and \mythetag{3.9} again mean that the 
cubic terms in \mythetag{2.2} are proportional to each other. Therefore, 
applying the left composition with some properly chosen linear transformation 
of the form \mythetag{2.1}, we can pass to an equivalent cubic transformation 
such that $F^{\kern 1pt 2}_{111}=F^{\kern 1pt 2}_{112}=F^{\kern 1pt 2}_{122}
=F^{\kern 1pt 2}_{222}=0$.\par
     If the inequality in \mythetag{3.8} is not fulfilled, then the first two
columns \mythetag{3.3} are zero. In this case at least one of the last two 
columns \mythetag{3.3} is nonzero. Assume that 
$$
\xalignat 3
&\hskip -2em
\hskip -2em
\Vmatrix F^{\kern 1pt 1}_{111}\\ 
\vspace{2ex}F^{\kern 1pt 2}_{111}\endVmatrix=0,
&&\Vmatrix F^{\kern 1pt 1}_{112}\\ 
\vspace{2ex} F^{\kern 1pt 2}_{112}\endVmatrix=0,
&&\Vmatrix F^{\kern 1pt 1}_{122}\\ 
\vspace{2ex} F^{\kern 1pt 2}_{122}\endVmatrix\neq 0.
\mytag{3.10}
\endxalignat
$$
Due to \mythetag{3.10} from $G_{2222}=0$ in \mythetag{3.1} we derive
$$
\hskip -2em
\Vmatrix F^{\kern 1pt 1}_{222}\\ 
\vspace{2ex} F^{\kern 1pt 2}_{222}\endVmatrix
=\alpha\Vmatrix F^{\kern 1pt 1}_{122}\\ 
\vspace{2ex}F^{\kern 1pt 2}_{122}\endVmatrix. 
\mytag{3.11}
$$
The equality \mythetag{3.11} along with \mythetag{3.10} means that  the 
cubic terms in \mythetag{2.2} are proportional to each other. Therefore, 
applying the left composition with some properly chosen linear transformation 
of the form \mythetag{2.1}, we can pass to an equivalent cubic transformation 
such that $F^{\kern 1pt 2}_{111}=F^{\kern 1pt 2}_{112}=F^{\kern 1pt 2}_{122}
=F^{\kern 1pt 2}_{222}=0$.\par
     And finally, if the inequality in \mythetag{3.10} is not fulfilled, then
the first three columns in \mythetag{3.3} are zero, while the last column is
nonzero. In this case applying the left composition with some properly chosen 
linear transformation of the form \mythetag{2.1}, we can pass to an equivalent 
cubic transformation such that 
$$
\hskip -2em
F^{\kern 1pt 2}_{111}=F^{\kern 1pt 2}_{112}
=F^{\kern 1pt 2}_{122}=F^{\kern 1pt 2}_{222}=0.
\mytag{3.12}
$$
The above results are summarized in the following theorem.
\mytheorem{3.1} In the case where the associated quartic form	\mythetag{2.4} 
is zero any cubic transformation \mythetag{1.1} is equivalent to a cubic
transformation whose leading coefficients in \mythetag{2.2} obey the
relationships \mythetag{3.12}.
\endproclaim 
\head
4. The case where $\omega[1]$ is indefinite.
\endhead
     An indefinite form takes both positive and negative values. Therefore
it takes zero values either. The form $\omega[1]$ is a quartic form in
$\Bbb R^2$. It is given by the formula \mythetag{2.4}. Without loss of
generality we can assume that $G_{1111}\neq 0$. Indeed, otherwise we can 
apply Theorem~\mythetheorem{2.1} and pass to an equivalent cubic 
transformation whose determinant $G_{1111}\neq 0$. The property of 
the associated quartic form $\omega[1]$ being indefinite is invariant 
under passing to an equivalent cubic transformation. This fact follows from 
Theorems~4.1 and 4.2 in \mycite{12}.\par 
     Thus we have an indefinite quartic form $\omega[1]$ in \mythetag{2.4}
with $G_{1111}\neq 0$. Setting $z^2=1$, we can treat $\omega[1]$ as a 
quartic polynomial with respect to $z^1$. This polynomial takes both 
positive and negative values. It is an elementary fact that a polynomial
of even degree taking values of both signs has at least two different 
roots $z^1\neq\tilde z^1$. Due to $z^1\neq\tilde z^1$ two vector-columns
$$
\xalignat 2
&\hskip -2em
\Vmatrix z^1\\ 
\vspace{2ex} 1\endVmatrix,
&&\Vmatrix \tilde z^1\\
\vspace{2ex} 1\endVmatrix
\mytag{4.1}
\endxalignat
$$
are linearly independent. Applying Theorem~\mythetheorem{2.1} and the 
formulas \mythetag{2.10} to the vector-columns \mythetag{4.1}, we define 
the matrix $T$ by setting
$$
\xalignat 2
&\hskip -2em
T^1_1=z^1,
&&T^1_2=\tilde z^1,\\
\vspace{-1.5ex}
\mytag{4.2}\\
\vspace{-1.5ex}
&\hskip -2em
T^2_1=1,
&&T^2_2=1.
\endxalignat
$$
The matrix $T$ defined through \mythetag{4.2} is non-degenerate. Using it
as the matrix of a linear transformation in \mythetag{2.1}, we can pass to 
an equivalent cubic transformation such that two of the six determinants
in \mythetag{2.3} do vanish simultaneously
$$
\xalignat 2
&\hskip -2em
G_{1111}=0,
&&G_{2222}=0. 
\mytag{4.3}
\endxalignat
$$
Along with \mythetag{4.3}, some of the other four determinants can 
vanish. Vanishing of them specify refinement subcases within the case 
of an indefinite form $\omega[1]$. 
\mydefinition{4.1} A cubic transformation $f$ of $\Bbb R^2$ with indefinite
form $\omega[1]$ is called obeying the first refinement condition (or the
condition R1) if there is an equivalent cubic transformation $\tilde f$
such that the equalities \mythetag{4.3} are fulfilled along with 
$$
\hskip -2em
G_{1112}=0.
\mytag{4.4}
$$
\enddefinition
\mydefinition{4.2} A cubic transformation $f$ of $\Bbb R^2$ with indefinite
form $\omega[1]$ is called obeying the second refinement condition (or the
condition R2) if there is an equivalent cubic transformation $\tilde f$ such 
that the equalities \mythetag{4.3} are fulfilled along with 
$$
\hskip -2em
G_{1122}=0.
\mytag{4.5}
$$
\enddefinition
\mydefinition{4.3} A cubic transformation $f$ of $\Bbb R^2$ with indefinite
form $\omega[1]$ is called obeying the third refinement condition (or the
condition R3) if there is an equivalent cubic transformation $\tilde f$ such 
that the equalities \mythetag{4.3} are fulfilled along with 
$$
\pagebreak
\hskip -2em
G_{1212}=0.
\mytag{4.6}
$$
\enddefinition
\mydefinition{4.4} A cubic transformation $f$ of $\Bbb R^2$ with indefinite
form $\omega[1]$ is called obeying the fourth refinement condition (or the
condition R4) if there is an equivalent cubic transformation $\tilde f$ such 
that the equalities \mythetag{4.3} are fulfilled along with 
$$
\hskip -2em
G_{1222}=0.
\mytag{4.7}
$$
\enddefinition
The equalities \mythetag{4.4}, \mythetag{4.5}, \mythetag{4.6}, \mythetag{4.7}
can be fulfilled either separately or simultaneously in various combinations.
Therefore one can speak of the conditions R1.2, R1.3, R1.4, R1.2.3 etc. 
Studying all of them is far beyond the scope of the present paper. We note
only that the condition R1.2.3.4 is incompatible with the case of an indefinite 
form $\omega[1]$ since it implies $\omega[1]=0$.\par
\head
5. The case where $\omega[1]$ is semi-definite.
\endhead
     A semi-definite form takes values of some definite sign and vanishes for
some values of its arguments not all zero. The property of $\omega[1]$ being
semi-definite is invariant under passing to an equivalent cubic transformation. 
Therefore, like in the previous section, without loss of generality we can 
initially assume that $G_{1111}\neq 0$. Then we set $z^2=1$ and treat 
$\omega[1]$ as a quartic polynomial with respect to $z^1$. This polynomial 
takes values of the same sign and has real roots since $\omega[1]$ is 
semi-definite. The nature of its roots defines two subcases of the present 
case:
\roster
\item"1)" the case where $\omega[1](z^1,1)$ has two distinct roots each of 
multiplicity two;
\item"2)" the case where $\omega[1](z^1,1)$ has one root of multiplicity four. 
\endroster
Actually the subcases 1) and 2) can be defined in other words. Indeed, 
$\omega[1](z^1,z^2)$ in \mythetag{2.4} is a homogeneous quartic polynomial 
of two variables $z^1$ and $z^2$. It can be understood as a function in the
real projective space $\RP^1$. Then the above two subcases are formulated as 
follows: 
\roster
\item"1)" the case where $\omega[1]$ has two distinct roots in $\RP^1$ each of 
multiplicity two;
\item"2)" the case where $\omega[1]$ has one root of multiplicity four in $\RP^1$. 
\endroster\par
     The subcase 1) of the present case is similar to the case considered
in the previous section. Indeed, if $(z^1,z^2)$ and $(\tilde z^1,\tilde z^2)$
are two distinct roots of the form $\omega[1]$ in $\RP^1$, then we can use them
for defining a non-degenerate matrix $T$: 
$$
\xalignat 2
&\hskip -2em
T^1_1=z^1,
&&T^1_2=\tilde z^1,\\
\vspace{-1.5ex}
\mytag{5.1}\\
\vspace{-1.5ex}
&\hskip -2em
T^2_1=z^2,
&&T^2_2=\tilde z^2.
\endxalignat
$$
Applying the matrix $T$ with the components \mythetag{5.1} as the matrix of a
linear transformation in \mythetag{2.1}, we can pass to an equivalent cubic
transformation such that 
$$
\xalignat 2
&\hskip -2em
G_{1111}=0 
&&G_{2222}=0. 
\mytag{5.2}
\endxalignat
$$
The equalities \mythetag{5.2} do coincide with \mythetag{4.3}. Therefore we can
complement them with the refinement conditions R1, R2, R3, R4 using definitions
similar to Defini\-tions~\mythedefinition{4.1}, \mythedefinition{4.2},  
\mythedefinition{4.3}, and ~\mythedefinition{4.4}. These refinement conditions
can be combined into groups forming the conditions R1.2, R1.3, R1.4, R1.2.3 etc
as described above. The condition R1.2.3.4 is incompatible with the present 
case of a semi-definite\linebreak form $\omega[1]$ since it implies \pagebreak 
$\omega[1]=0$.\par
     The subcase 2) is more complicated. In this subcase the root $(z^1,z^2)$
of the form $\omega[1]$ in $\RP^1$ defines only the first column of the 
matrix $T$ in \mythetag{5.1}. The second column is deliberate or possibly it can 
be specified using the other forms $\omega[2]$, $\omega[3]$, $\omega[4]$,
$\omega[5]$. More details are to be elaborated in forthcoming papers.\par
\head
5. The case where $\omega[1]$ is definite.
\endhead
     A definite form takes values of some definite sign and does not vanish
if its arguments are not all zero. The property of the form $\omega[1]$ being
definite is invariant under passing to an equivalent cubic transformation. 
This fact follows from Theorems~4.1 and 4.2 in \mycite{12}. Other details of
the case of a definite form $\omega[1]$ are yet uncertain. They should be 
elaborated in future within a refinement procedure.\par
\head
5. Conclusions.
\endhead
     The rough classification scheme suggested in the present paper is only 
a framework that does not contain all possible subcases and canonical presentations
of cubic transformations within each particular subcase. These subcases and canonical
presentations should be discovered and cataloged in forthcoming papers within a 
refinement procedure.\par

\enddocument
\end